\documentclass{article}

\usepackage[caption=false]{subfig}
\usepackage{amsmath,amssymb,amsfonts,amsthm,caption}
\usepackage{graphics}
\usepackage{epstopdf}
\usepackage{bbm}
\usepackage{enumerate}
\usepackage{epsfig}
\usepackage{booktabs,multirow}
\usepackage{url,xcolor}

\usepackage{a4wide}

\theoremstyle{plain}

\theoremstyle{definition}

\theoremstyle{remark}
\newtheorem{remark}{Remark}

\renewcommand{\phi}{\varphi}

\newcommand{\beqn}{\begin{eqnarray}}
\newcommand{\eeqn}{\end{eqnarray}}
\newcommand{\beqnn}{\begin{eqnarray*}}
\newcommand{\eeqnn}{\end{eqnarray*}}

\newcommand{\Dl}{{}^{\scriptscriptstyle C}_{\scriptscriptstyle 0}\!D^{\alpha}_{\scriptscriptstyle t}}
\newcommand{\DlRL}{{}_{\scriptscriptstyle 0}D^{\alpha}_{\scriptscriptstyle t'}}
\newcommand{\DrRL}{{}_{\scriptscriptstyle t}D^{\alpha}_{\scriptscriptstyle t_f}}

\usepackage{dsfont}

\begin{document}

\title{Optimal control of a fractional order epidemic model
with application to human respiratory syncytial virus infection\thanks{This is a preprint 
of a paper whose final and definite form is with \emph{Chaos, Solitons \& Fractals}, 
available from {\tt http://www.elsevier.com/locate/issn/09600779}. 
Submitted 23-July-2018; Revised 14-Oct-2018; Accepted 15-Oct-2018.}}

\author{Silv\'erio Rosa$^a$\\
\texttt{rosa@ubi.pt}
\and
Delfim F. M. Torres$^b$\thanks{Corresponding author.}\\
\texttt{delfim@ua.pt}}

\date{$^a$Instituto de Telecomunica\c{c}\~{o}es
and Department of Mathematics,\\
Universidade da Beira Interior,
6201-001 Covilh\~a, Portugal\\[0.3cm]
$^b$Center for Research and Development in Mathematics and Applications (CIDMA),
Department of Mathematics,\\
University of Aveiro, 3810-193 Aveiro, Portugal}

\maketitle


\begin{abstract}
A human respiratory syncytial virus surveillance system was implemented
in Florida in 1999, to support clinical decision-making for prophylaxis
of premature newborns. Recently, a local periodic SEIRS mathematical model was
proposed in [Stat. Optim. Inf. Comput. 6 (2018), no.~1, 139--149]
to describe real data collected by Florida's system.
In contrast, here we propose a non-local fractional (non-integer) order model.
A fractional optimal control problem is then formulated and solved,
having treatment as the control. Finally, a cost-effectiveness analysis
is carried out to evaluate the cost and the effectiveness of
proposed control measures during the intervention period,
showing the superiority of obtained results with respect to previous ones.

\medskip


\medskip

\end{abstract}


\section{Introduction}

Human respiratory syncytial virus (HRSV) is a virus that causes
respiratory tract infections. It is a major cause of lower respiratory
tract infections and hospital visits during infancy and childhood.
A prophylactic medication, palivizumab, can be employed to prevent HRSV in pre-term
infants (under 35 weeks gestation), infants with certain congenital heart defects
or bronchopulmonary dysplasia, and infants with congenital malformations
of the airway. Treatment is limited to supportive care, including oxygen therapy.
In temperate climates, there is an annual epidemic during the winter season.
In tropical climates, infection is most common during the rainy season.
In the United States, 60\% of infants are infected during their first HRSV season,
and nearly all children will have been infected with the virus by two to three years of age \cite{Glezen}.
Of those infected with HRSV, 2 to 3\% will develop bronchiolitis, necessitating hospitalization \cite{Hall}.
Natural infection with HRSV induces protective immunity, which wanes over time, possibly
more so than other respiratory viral infections, and thus people can be infected multiple times.
Sometimes, an infant can become symptomatically infected more than once, even within a single HRSV season.
Severe HRSV infections have increasingly been found among elderly patients.
Young adults can be re-infected every five to seven years, with symptoms looking
like a sinus infection or a cold. The Florida Department of Health provides an integrated
and reliable HRSV system, with data from hospitals and laboratories \cite{flhealth}.

Mathematical models can project how infectious diseases progress,
to show the likely outcome of an epidemic, and help inform public health
interventions. In epidemiology, compartmental models serve as the base
mathematical framework for understanding the complex dynamics of these systems.
Such compartments, in the simplest case, stratify the population into two health states:
susceptible to the infection of the pathogen, often denoted by $S$, and infected by the
pathogen, often denoted by the symbol $I$. The way that these compartments interact
is based upon phenomenological assumptions, and the model is built up from there.
These models are usually investigated through ordinary differential equations.
Depending on the disease, other compartments may be included,
most notably the recovered/removed/immune compartment, often denoted by $R$. 
Recently, to push such models to further realism and taking into account
the influence of past on the current and future state of the diseases,
they have been characterized mathematically with the help of
fractional order differential equations: see, e.g.,
investigations in dengue \cite{doi:10.1063/1.3636838},
Ebola \cite{area2015mathematical}, tuberculosis \cite{TB:frac:2018} and
HIV/AIDS \cite{MR3743014,MR3691319}.

A crucial question consists to find parameters
for the particular disease under study, and use those parameters to calculate
the effects of possible control interventions, like treatment or vaccination.
Then the central issue is how to implement such interventions in an optimal way.
This investigation program has been carried out for several infectious diseases
with classical integer-order compartmental models: see, e.g., \cite{rosa:delfim2018parameter}
for HRSV and \cite{ndairou2017mathematical} and \cite{MR3810766,MR3731603,MR3703345} 
for Zika and Ebola viruses, respectively, where the implementation 
of optimal control interventions are proposed.
Here we extend such approach with new fractional compartmental models.
Applications of fractional calculus in numerous fields of science and
engineering have gained popularity and importance in past four decades, 
see \cite{agarwal2015extended,ruzhansky2017advances,Saoudi2018} 
and references therein. Recently, extensions of fractional derivative 
operators have been developed and proved to be very useful 
in several applications, showing a high vitality 
of the research field \cite{Agarwal2017,agarwal2017extended,MR3797763}.
An extension of the Caputo fractional derivative operator
is given in \cite{MR3756437} by using a generalized beta function,
Saigo--Maeda fractional differential operators involving Appell functions 
are investigated in \cite{MR3764705}, while fractional integral operators 
involving Gauss hypergeometric functions are studied in \cite{choi2015certain}.
For variable-order fractional operators see \cite{MR3822307}
and references therein. Fractional operators on arbitrary time scales 
are proposed and investigated in \cite{MR2800417,MR3571716}.
Usually, numerical techniques for solving such fractional order models
are required \cite{AGARWAL201840,MR3787702}. To the best of our knowledge, 
this is the first work to use fractional calculus and fractional 
optimal control in the study of HRSV. For that we use derivatives
in the standard Caputo sense.

A comparison of the standard SIRS model
with a more complex integer-order model of HRSV transmission,
in which individuals acquire immunity gradually
after repeated exposure to infection, is given in \cite{Weber2001}.
In \cite{MR2718412}, an age-structured mathematical model for HRSV
is proposed, where children younger than one year old,
who are the most affected by this illness, are specially considered.
A numerical scheme for the SIRS seasonal epidemiological
model of HRSV transmission is proposed in \cite{MR2435573}. It turns out
that solutions for HRSV compartmental models are typically memory-periodic
systems \cite{MR2426325}. For this reason, in this work we propose
the use of fractional optimal control theory \cite{MR3673702,MR3673710} to a
non-autonomous fractional SEIRS model, and show its usefulness
according with real HRSV data provided by the Florida Department
of Health \cite{flhealth}.

The paper is organized as follows. In Section~\ref{sec:2},
we introduce the SIRS-$\alpha$ and SEIRS-$\alpha$ fractional epidemic models,
which generalize corresponding integer-order ($\alpha = 1$) models of \cite{UBI:UA}.
The main results are then given in Section~\ref{sec:3}: estimation of the fractional
order $\alpha$ with real data of Florida, for the two proposed fractional models
(Section~\ref{sec:alphaorder}); fractional optimal control, cost-effectiveness
and numerical simulations for the more realistic SEIRS-$\alpha$ model
(Sections~\ref{sec:optctrl} and \ref{sec:numresul}).
We end with Section~\ref{sec:conclu} of conclusions
and perspectives of future work.


\section{Fractional compartmental models}
\label{sec:2}

We focus on compartmental models that divide the population into mutually
exclusive distinct groups (of susceptible, or infected, or immune individuals,
or \ldots) and we use  deterministic continuous transitions between those groups,
also known as states. Due to the seasonality of HRSV, the models that best fit
real data are periodic. In \cite{Weber2001},
two integer-order models are proposed, where the transmission is periodic:
(i) a simple model with only three compartments, known as SIRS, which we extend
here to the fractional SIRS-$\alpha$ case in Section~\ref{subsec:SIRS};
(ii) and a more complex model with seventeen compartments, named MSEIRS4. However,
it is shown that the simpler $SIRS$ model fits better real data \cite{Weber2001}.
Zang et al. \cite{zhang2012existence} use a non-autonomous SEIR model where,
beyond the periodicity in the transmission rate, the annual recruitment
rate is also periodic. This assumption is due to opening and closing
of schools. Here we generalize such ideas and consider
a simple periodic fractional SEIRS-$\alpha$ model (Section~\ref{subsec:SEIRS}).


\subsection{SIRS-$\alpha$ model}
\label{subsec:SIRS}

In our first fractional model, we consider that the population consists
of susceptible ($S$), infected and infectious ($I$), and recovered ($R$) individuals.
A characteristic feature of HRSV is that immunity after infection is temporary,
so that the recovered individuals become susceptible again \cite{Weber2001}.
Let parameter $\mu$ denote the birth rate, which we assume equal to the mortality rate;
$\gamma$ be the rate of loss of immunity; and $\nu$ the  rate of loss of infectiousness.
The influence of the seasonality on the transmission parameter $\beta$ is modeled by the
\emph{cosine} function. Using a linear mass fractional action law, we propose
the following system of fractional differential equations:
\begin{equation}
\label{eq:modSIRS}
\begin{cases}
\Dl S(t)=  \mu-\mu S(t)-\beta(t)S(t) I(t) +\gamma R(t),\\[1.2mm]
\Dl I(t)= \beta(t)S(t) I(t) -\nu I(t)-\mu I(t),\\[1.2mm]
\Dl R(t)=\nu I(t)-\mu R(t)-\gamma R(t),
\end{cases}
\end{equation}
where  $\beta(t)=b_0(1+b_1\cos(2 \pi t +\Phi))$ and $\Dl$
denotes the left Caputo derivative of order $\alpha \in (0,1]$
\cite{MR1658022,MR3736617}. The parameter $b_0$ is the mean
of the transmission parameter and $b_1$ is the amplitude of the seasonal fluctuation
in the transmission parameter, $\beta$.

\begin{remark}
When $\alpha=1$, the fractional compartmental model \eqref{eq:modSIRS}
represent the classical SIRS model investigated in \cite{UBI:UA}.
\end{remark}

We improve previous paper \cite{rosa:delfim2018parameter}, 
where two models were used: an SIRS and a SEIRS model. The first model, 
periodic in the transmission rate, is the simplest, however, as already mentioned, 
it proven to be better than a MSEIRS4 model, which is a more complex model. 
The second model, adds to the former periodicity in the recruitment rate, 
fitting better the real data \cite{rosa:delfim2018parameter}.
Here we investigate SIRS and SEIRS like models, denoted by SIRS-$\alpha$ \eqref{eq:modSIRS}
and SEIRS-$\alpha$ \eqref{eq:modSEIRS}, where $\alpha$ denotes the non-integer order 
of differentiation under consideration, that are better than the (integer-order) 
models investigated in \cite{rosa:delfim2018parameter}. This makes our new $\alpha$-models 
interesting from a mathematical modeling point of view.


\subsection{SEIRS-$\alpha$ model}
\label{subsec:SEIRS}

To incorporate more features of HRSV, we extend the previous fractional-order
model in the following way. First, we include a latency period
by introducing a group $E$ of individuals who have been infected
but are not yet infectious. These individuals become
infectious at a rate $\varepsilon$. We assume the latency period to be
equal to the time between infection and the first symptoms.
As in \cite{zhang2012existence}, we consider that the annual
recruitment rate is seasonal due to schools opening/closing periods.
Our system of fractional differential equations is now given by
\begin{equation}
\label{eq:modSEIRS}
\begin{cases}
\Dl S(t) =  \lambda(t)-\mu S(t)-\beta(t) S(t) I(t) +\gamma R(t),\\[1.2mm]
\Dl E(t) = \beta(t) S(t) I(t) -\mu E(t)-\varepsilon E(t),\\[1.2mm]
\Dl I(t) =\varepsilon E(t)-\mu I(t)-\nu I(t),\\[1.2mm]
\Dl R(t) =\nu I(t)-\mu R(t)-\gamma R(t),
\end{cases}
\end{equation}
where $\lambda(t)=\mu(1 + c_1  \cos( 2 \pi t +\Phi) )$ is the recruitment rate
(including newborns and immigrants), parameter $c_1$ is the amplitude
of the seasonal fluctuation in the recruitment parameter, $\lambda$,
while $\Phi$ is an angle that is chosen in agreement with real data,
and, as before, $\Dl$ denotes the left Caputo derivative of order
$\alpha \in (0,1]$. Note that in the particular case $\alpha = 1$
we obtain from \eqref{eq:modSEIRS} the SEIRS model of \cite{UBI:UA}.


\section{Main results}
\label{sec:3}

We begin by investigating how realistic the fractional models discussed in
Section~\ref{sec:2} are, with respect to HRSV and real data from
Florida \cite{flhealth}. For that, we borrow the values
for the parameters $\mu$, $\nu$, $\gamma$, $\varepsilon$, $b_0$,
$b_1$, $c_1$, and $\Phi$ from \cite{rosa:delfim2018parameter}
and do a proper estimation of the best fractional order $\alpha$
for models SIRS-$\alpha$ and SEIRS-$\alpha$.


\subsection{Estimation of the fractional order $\alpha$ for Caputo derivatives}
\label{sec:alphaorder}

\begin{figure}[!htb]
\centering
\includegraphics[scale=0.7]{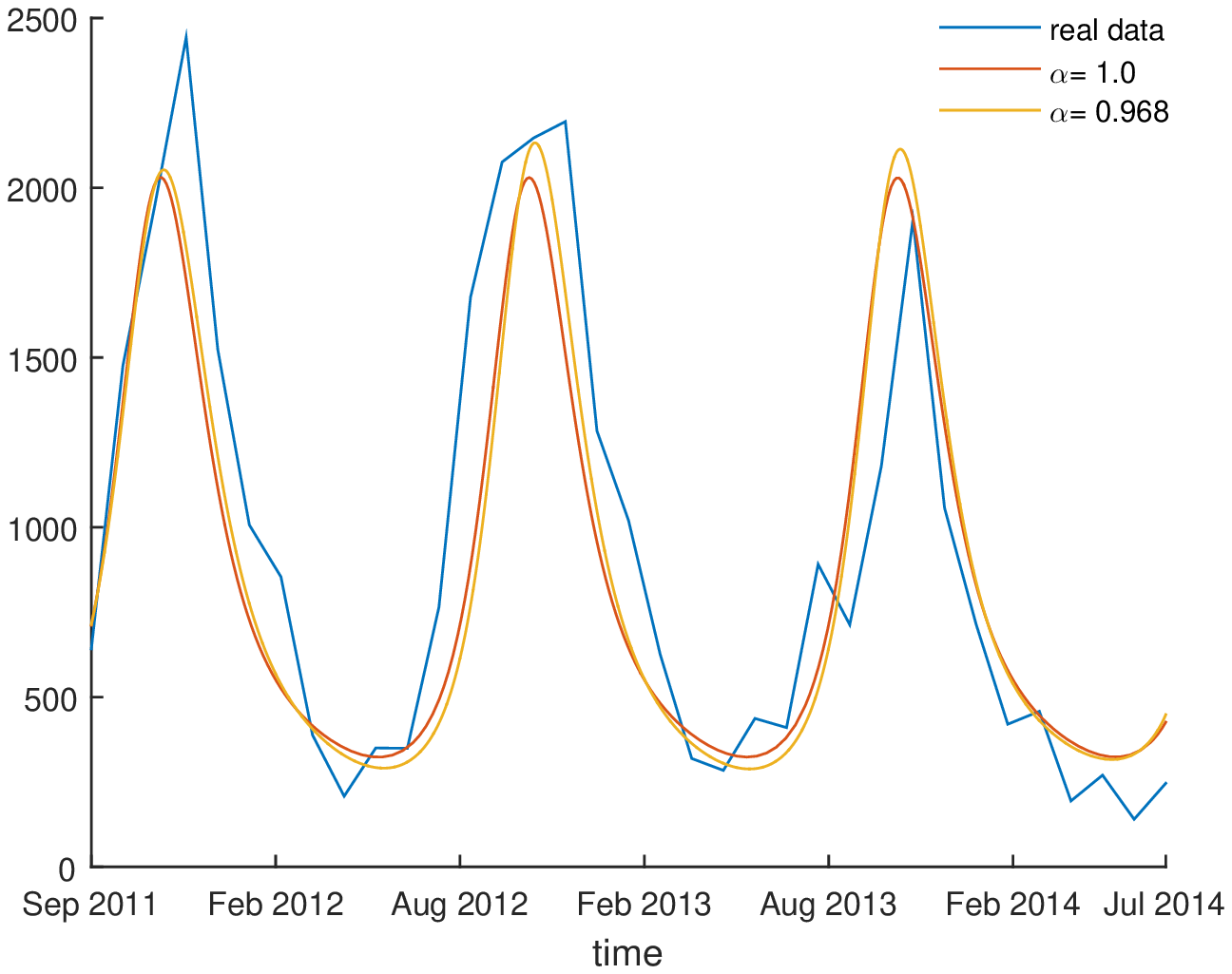}
\caption{Comparison of infected/infectious individuals $I(t)$:
real data, classical SIRS model (i.e., $\alpha = 1$)
and the fractional SIRS-$\alpha$ model with $\alpha=0.968$.}
\label{fig:sirs_alpha}
\end{figure}

\begin{figure}[!htb]
\centering
\includegraphics[scale=0.7]{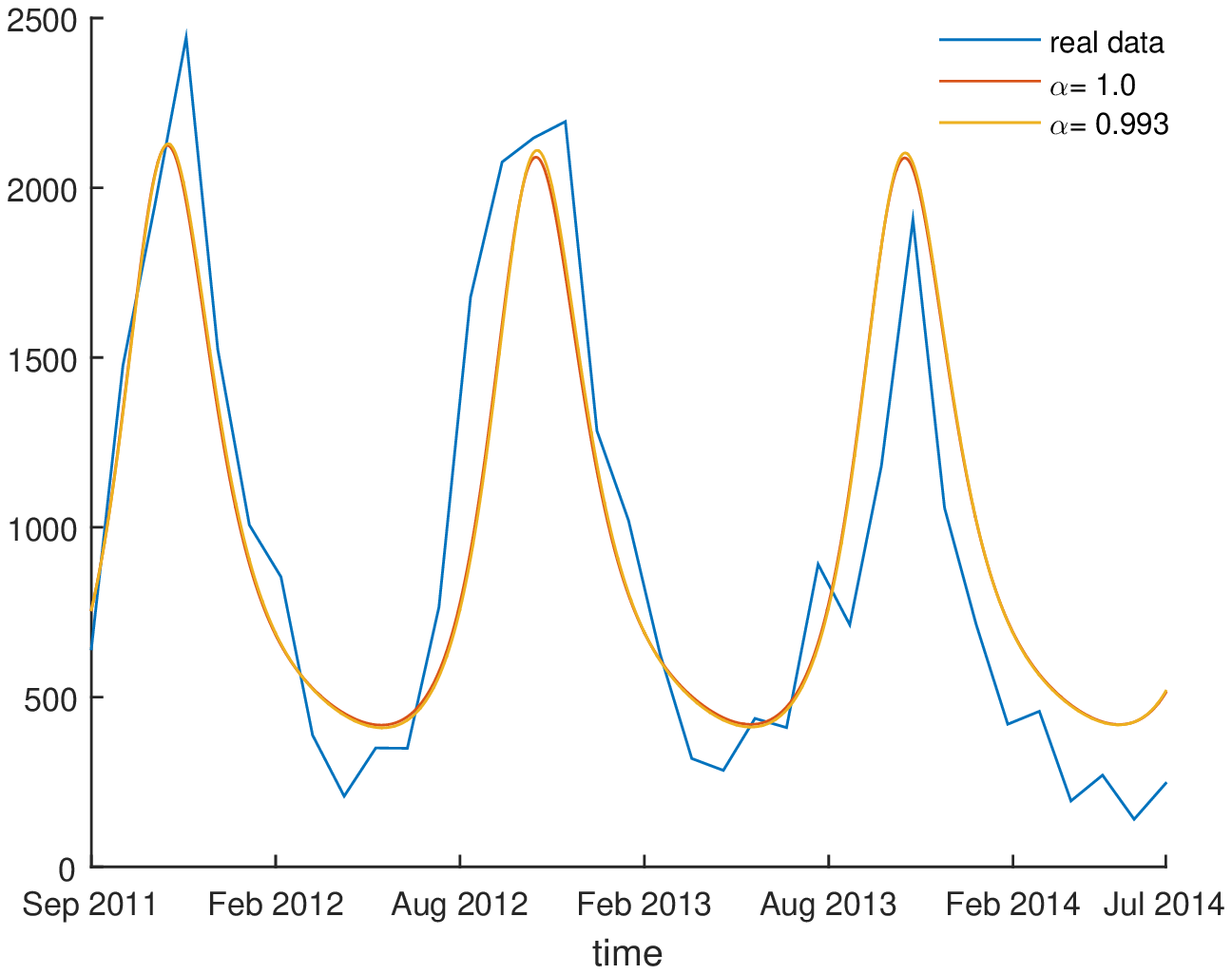}
\caption{Comparison of infected/infectious individuals $I(t)$:
real data, classical SEIRS model (i.e., $\alpha = 1$)
and the fractional SEIRS-$\alpha$ model with $\alpha=0.993$.}
\label{fig:seirs_alpha}
\end{figure}

Using the parameters values as given in Table~\ref{tab:param},
we searched the fractional order of differentiation, $\alpha$, that
best fits the data on the reported number of positive tests of
HRSV disease, per month, in the state of Florida, excluding North region,
between September 2011 and July 2014 (35 months). The data was obtained from
the Florida Department of Health \cite{flhealth}. 
The value of the derivative order, $\alpha$, was obtained by 
a search on the interval $]0,1]$. We started with $\alpha=1$ 
and successively lower its value until we find one whose 
lower values in the neighborhood correspond to worst 
fitting results. The results for the SIRS-$\alpha$ model 
\eqref{eq:modSIRS} are given in Figure~\ref{fig:sirs_alpha},
while the results corresponding to the SEIRS-$\alpha$ model 
\eqref{eq:modSEIRS} are given in Figure~\ref{fig:seirs_alpha}.
\begin{table}[!htb]
\caption{Models' parameters borrowed from \cite{rosa:delfim2018parameter}.}
\label{tab:param}
\centering
\begin{tabular}{ccccccccc}\toprule
model & $\mu$ & $\nu$ & $\gamma$ & $\varepsilon$ & $b_0$
& $b_1$ & $c_1$ & $\Phi$ \\[1mm] \midrule
SIRS-$\alpha$ & 0.0113 & 36 & 1.8 & --  & 74.2 & 0.14 & --  & $7\pi/5$ \\[1mm]
SEIRS-$\alpha$ & 0.0113 & 36 & 1.8 & 91  & 88.25 & 0.17 & 0.17 & $7\pi/5$ \\[1mm]
\bottomrule
\end{tabular}
\end{table}
\begin{table}[!htb]
\centering
\caption{Comparison of models \eqref{eq:modSIRS}
and \eqref{eq:modSEIRS} with real data.
The $l_2$ norm of the difference between real data
and the predictive cases given by the models is denoted by \texttt{error}.
In the last column, \texttt{relative error} denotes the percentual difference
of infants, per year, with respect to the total child population of Florida in 2014.}
\label{tab:errors}
\begin{tabular}{lllcccc}\toprule
model && \multicolumn{1}{c}{$\alpha$} && \texttt{error}
&& \texttt{relative error (\%)} \\[1mm] \midrule
\multirow{2}{*}{SIRS-$\alpha$}  && 1.0 && 1871.46 && 0.06567 \\[1mm]
&& 0.968    && 1840.90 && 0.06460    \\[1mm] \midrule
\multirow{2}{*}{SEIRS-$\alpha$} && 1.0      && 1719.12 && 0.06032\\[1mm]
&& 0.993    && 1716.91 && 0.06025 \\[1mm] \bottomrule
\end{tabular}
\end{table}

For comparison reasons, we adopted the same fitting approach
as in \cite{rosa:delfim2018parameter}, which consists in the
minimization of the $l_2$ norm of the difference between real data
and the number of HRSV infection individuals predicted by the models.
Table~\ref{tab:errors} gives the values of the errors for the considered
models. The relative errors show that these models fit quite well the data of HRSV disease.
Comparing the results for the optimal values of $\alpha$, that is,
$\alpha=0.968$ for the fractional order model \eqref{eq:modSIRS} and
$\alpha=0.993$ for the fractional order model \eqref{eq:modSEIRS},
with the classical SIRS and SEIRS models obtained for $\alpha = 1$,
one concludes that while the absolute error of the fractional SIRS-$\alpha$
model reduces more than the homonymous error of the fractional SEIRS-$\alpha$ model,
it is the SEIRS-$\alpha$ model the one with lowest error.


\subsection{Fractional optimal control problem}
\label{sec:optctrl}

In order to investigate some optimal measures, we choose,
in agreement with Section~\ref{sec:alphaorder},
the SEIRS-$\alpha$ model with $\alpha=0.993$, which is the one that
provides the best fitting to the considered real data. 
The evolution of the variables of the model depend on some 
circumstances that can be controlled. In what concerns
HRSV disease, treatment is the most commonly used. Hence, we consider the following fractional
optimal control problem: to minimize the number of infectious individuals and the cost
associated to control the disease with the treatment of the patients, that is,
\begin{equation}
\label{cost-functional}
\min ~\mathcal{J}(I,\mathbbm{T})
=\int_0^{t_f} \left(\kappa_1\,I(t)+ \kappa_2\, \mathbbm{T}^2(t)\right) ~dt
\end{equation}
with given $0<\kappa_1,\kappa_2 <\infty$, subject to the fractional control system
\begin{equation}
\label{eq:modSEIRS_control}
\begin{cases}
\Dl S(t) =  \lambda(t)-\mu S(t)-\beta(t) S(t) I(t) +\gamma R(t),\\
\Dl E(t) = \beta(t) S(t) I(t) -\mu E(t)-\varepsilon E(t),\\
\Dl I(t) =\varepsilon E(t)-\mu I(t)-\nu I(t)-\mathbbm{T}(t) I(t),\\
\Dl R(t) =\nu I(t)-\mu R(t)-\gamma R(t)+\mathbbm{T}(t) I(t)
\end{cases}
\end{equation}
and given initial conditions
\begin{equation}
\label{ocp:ic}
S(0),E(0),I(0),R(0)\geqslant 0.
\end{equation}
Here, $\mathbbm{T}$ is the control variable, which designates \emph{treatment}.
Note that in absence of treatment, that is, for $\mathbbm{T}(t) \equiv 0$,
then the control system \eqref{eq:modSEIRS_control} reduces to the SEIRS-$\alpha$
dynamical system \eqref{eq:modSEIRS}. The set of admissible control functions is
\begin{equation}
\label{Omega:set}
\Omega=\left\{\mathbbm{T}(\cdot)\in L^{\infty}(0,t_f):
0\leqslant \mathbbm{T}(t)\leqslant \mathbbm{T}_{\max},\forall t\in[0,t_f]\right\}.
\end{equation}
Pontryagin's maximum principle (PMP) for fractional optimal control can be used
to solve the problem \cite{MR3529374,MR3443073,MR3225198,malinowska2012introduction}.
The Hamiltonian of our optimal control problem is
\begin{multline*}
\mathcal{H}
= \kappa_1 I +\kappa_2 \mathbbm{T}^2+p_1(\lambda-\mu S-\beta S I +\gamma R)
+p_2(\beta S I -\mu E-\varepsilon E)\\
+p_3(\varepsilon E-\mu I-\nu I-\mathbbm{T}I)
+p_4(\nu I-\mu R-\gamma R+\mathbbm{T}I);
\end{multline*}
the optimality condition of PMP ensures that the optimal control is given by
\begin{equation}
\label{eq:ext:cont}
\mathbbm{T}(t)=\min\left\{\max\left\{0,\dfrac{(p_3(t)-p_4(t)) I(t)}{2
\kappa_2}\right\},\mathbbm{T}_{\max}\right\};
\end{equation}
while the adjoint system asserts that the co-state variables
$p_i(t)$, $i = 1,\ldots, 4$, satisfy
\begin{align}
\label{eq:co_states_system}
\begin{cases}
\DrRL \, p_1(t) = p_1(t) \left(\mu+\beta(t)I(t)\right)-\beta(t) I(t) p_2(t),\\[2mm]
\DrRL \, p_2(t) =  p_2(t) \left(\mu+\varepsilon\right)-\varepsilon p_3(t), \\[2mm]
\DrRL \, p_3(t)  =  -\kappa_1+\beta(t) p_1(t) S(t)-p_2(t) \beta(t) S(t)\\
\qquad\qquad\quad +p_3(t) \left(\mu+\nu+\mathbbm{T}(t)\right)
-p_4(t)\left(\nu+\mathbbm{T}(t)\right),\\[2mm]
\DrRL \, p_4(t) = -\gamma p_1(t)+p_4(t)\left(\mu+\gamma\right),
\end{cases}
\end{align}
which is a fractional system of right Riemann--Liouville derivatives $\DrRL$.
In addition, the following transversality conditions hold:
\begin{equation}
\label{eq:transversality}
_{\scriptscriptstyle t}D^{\alpha-1}_{\scriptscriptstyle t_f}
p_i\bigm|_{\scriptscriptstyle t_f}=0
\Leftrightarrow _{\scriptscriptstyle t}\!I^{1-\alpha}_{\scriptscriptstyle t_f}
p_i\bigm|_{\scriptscriptstyle t_f}=p_i(t_f)=0, \quad i=1,\ldots,4,
\end{equation}
where $_{\scriptscriptstyle t}\!I^{1-\alpha}_{\scriptscriptstyle t_f}$
is the right Riemann--Liouville fractional integral of order $1-\alpha$.


\subsection{Numerical results and cost-effectiveness analysis
of the fractional HRSV optimal control problem}
\label{sec:numresul}

The optimal control problem \eqref{cost-functional}--\eqref{Omega:set}
is numerically solved with the help of PMP and its optimality conditions,
as discussed in Section~\ref{sec:optctrl}, in the classical ($\alpha = 1$)
and fractional ($\alpha=0.993$) cases, using the
predict-evaluate-correct-evaluate (PECE) method of Adams--Basforth--Moulton
\cite{diethelm2005algorithms}. First we solve system \eqref{eq:modSEIRS_control}
by the PECE procedure with initial conditions for the state variables \eqref{ocp:ic}
given in terms of percentage of total population, that is, $S(0) + E(0) + I(0) + R(0) = 1$,
and a guess for the control over the time interval $[0,t_f]$,
and obtain the values of the state variables $S$, $E$, $I$ and $R$.
Applying the change of variable
\[
t'=t_f-t
\]
to the system of adjoint equations \eqref{eq:co_states_system}
and to the transversality conditions \eqref{eq:transversality},
we obtain the following left Riemann--Liouville fractional initial value problem
\eqref{eq:co_states_fr2}--\eqref{eq:trans2}:
\begin{equation}
\label{eq:co_states_fr2}
\begin{cases}
\DlRL\,  p_1(t') = -\left[p_1(t')(\mu+\beta(t')I(t'))-\beta(t') I(t') p_2(t')\right],\\[2mm]
\DlRL\, p_2(t')=  -\left[p_2(t')(\mu+\varepsilon)-\varepsilon p_3(t')\right], \\[2mm]
\DlRL\, p_3(t')  =  -\left[-\kappa_1+\beta(t') p_1(t') S(t')
+p_3(t')(\mu+\nu+\mathbbm{T}(t'))- p_4(t')(\nu+\mathbbm{T}(t'))\right],\\[2mm]
\DlRL\, p_4(t') = -\left[-\gamma p_1(t')+p_4(t')(\mu+\gamma)\right]
\end{cases}
\end{equation}
with initial conditions
\begin{equation}
\label{eq:trans2}
p_i(t')\bigm|_{\scriptscriptstyle t'=0}=0, \quad i=1,\ldots,4.
\end{equation}
Given the initial conditions \eqref{eq:trans2},
we solve \eqref{eq:co_states_fr2} with the PECE procedure and
obtain the values of the co-state variables $p_i$, $i=1,\ldots,4$.
The control is then updated by a convex combination of the previous control
and the value from \eqref{eq:ext:cont}. This procedure is repeated
iteratively until the values of the controls at the previous iteration
are very close to the ones at the current iteration.
To validate this algorithm, a fractional optimal control problem whose exact solution
is known, Example~3.1 in \cite[p.~86]{malinowska2012introduction},
was first successfully solved with it. Here we present our numerical
results to the HRSV optimal control problem \eqref{cost-functional}--\eqref{Omega:set},
for which an analytical solution is unknown.
\begin{table}[!htb]
\centering
\caption{Initial conditions \eqref{ocp:ic}, in terms of percentage of total population,
for the fractional optimal control problem \eqref{cost-functional}--\eqref{Omega:set}
with parameters given by Table~\ref{tab:param}, excepting
angle $\Phi$ that is here assumed to be $\pi/2$ so
that the values correspond to the endemic equilibrium
of \eqref{eq:modSEIRS}.}\label{tab:solinit}
\begin{tabular}{c@{\hspace*{2.2cm}}c@{\hspace*{2.2cm}}c@{\hspace*{2.2cm}}c}\toprule
$S(0)$ & $E(0)$ & $I(0)$ & $R(0)$ \\[1mm] \midrule
$0,4081$ & $0,0110$ &  $0,0278$ &  $0,5531$ \\ \bottomrule
\end{tabular}
\end{table}

In our numerical computations, we consider that $\mathbbm{T}_{\max}=1$
and the other parameters are fixed according to Table~\ref{tab:param},
with exception of angle $\Phi$ that is assumed to be $\pi/2$.
Such value allows the transmission parameter initial value to be the
average, $\beta(0)=b_0$, and the recruitment rate initial value
to be also the average, $\lambda(0)=\mu$. Our initial conditions,
given by Table~\ref{tab:solinit}, guarantee the existence of a
non-trivial endemic equilibrium for the system \eqref{eq:modSEIRS_control} without
control ($\mathbbm{T}(t) \equiv 0$), corresponding to the population system
\eqref{eq:modSEIRS} prior introduction of treatment.
Because World Health Organization goals for most diseases are usually
fixed for five years periods, we assumed $t_f=5$.
\begin{figure}[!htb]
\centering
\subfloat[Variation of the no. of susceptible individuals.]{\includegraphics[scale=0.46]{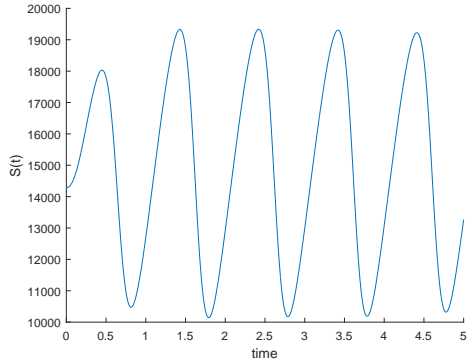}}
\hspace*{1cm}
\subfloat[Variation of the no. of exposed individuals.]{\includegraphics[scale=0.46]{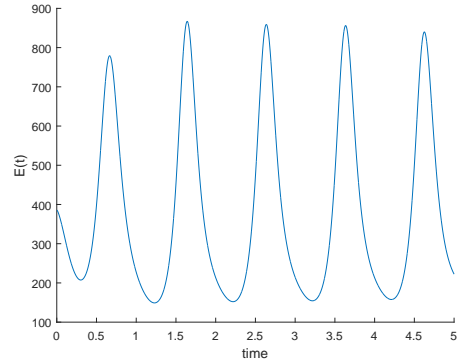}}\\
\subfloat[Variation of the no. of infected individuals.]{\includegraphics[scale=0.46]{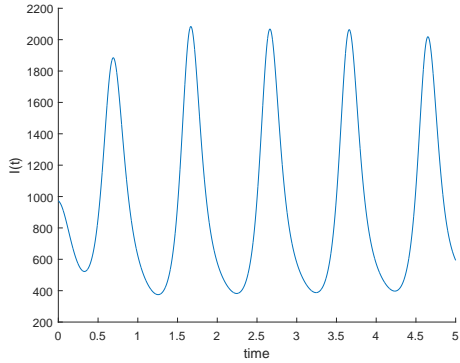}}
\hspace*{1cm}
\subfloat[Variation of the no. of recovered individuals.]{\includegraphics[scale=0.46]{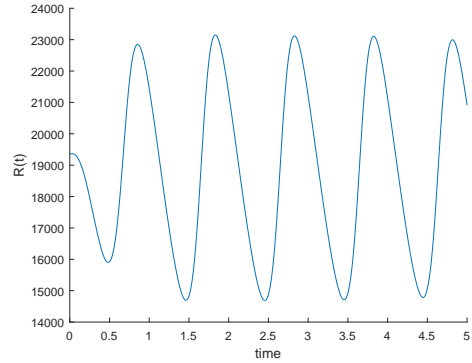}}
\caption{State variables of the fractional optimal control problem \eqref{cost-functional}--\eqref{Omega:set},
considering  $\alpha=0.993$ and weights $\kappa_1=1$ and $\kappa_2=0.001$.}
\label{fig:states_var:K}
\end{figure}
\begin{figure}[!htb]
\centering
\subfloat[Evolution of the four co-state variables.]{
\includegraphics[scale=0.46]{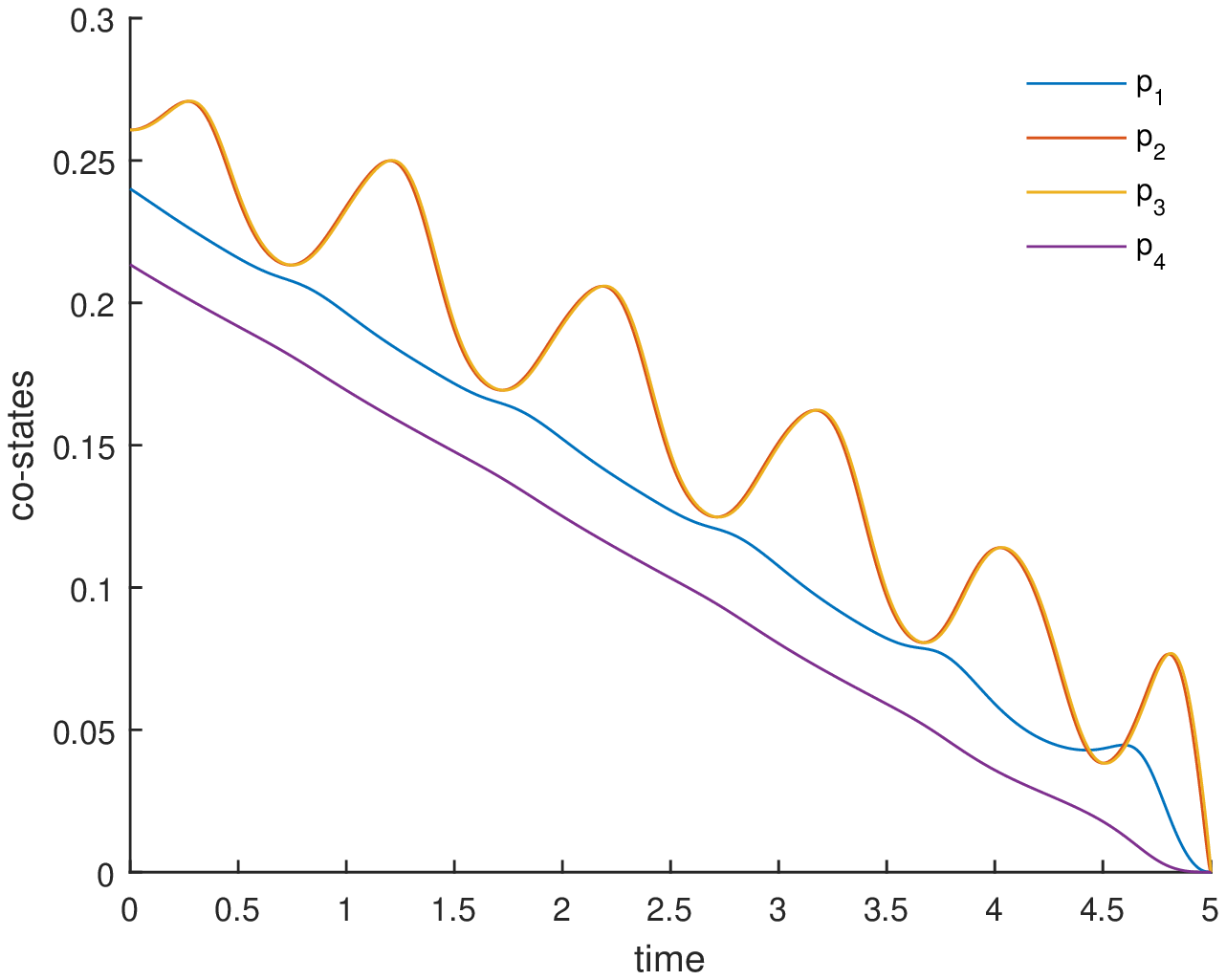}\label{fig:co_states_001}}
\hspace*{1cm}
\subfloat[Evolution of the optimal control $\mathbbm{T}$ (treatment).]{\includegraphics[scale=0.46]{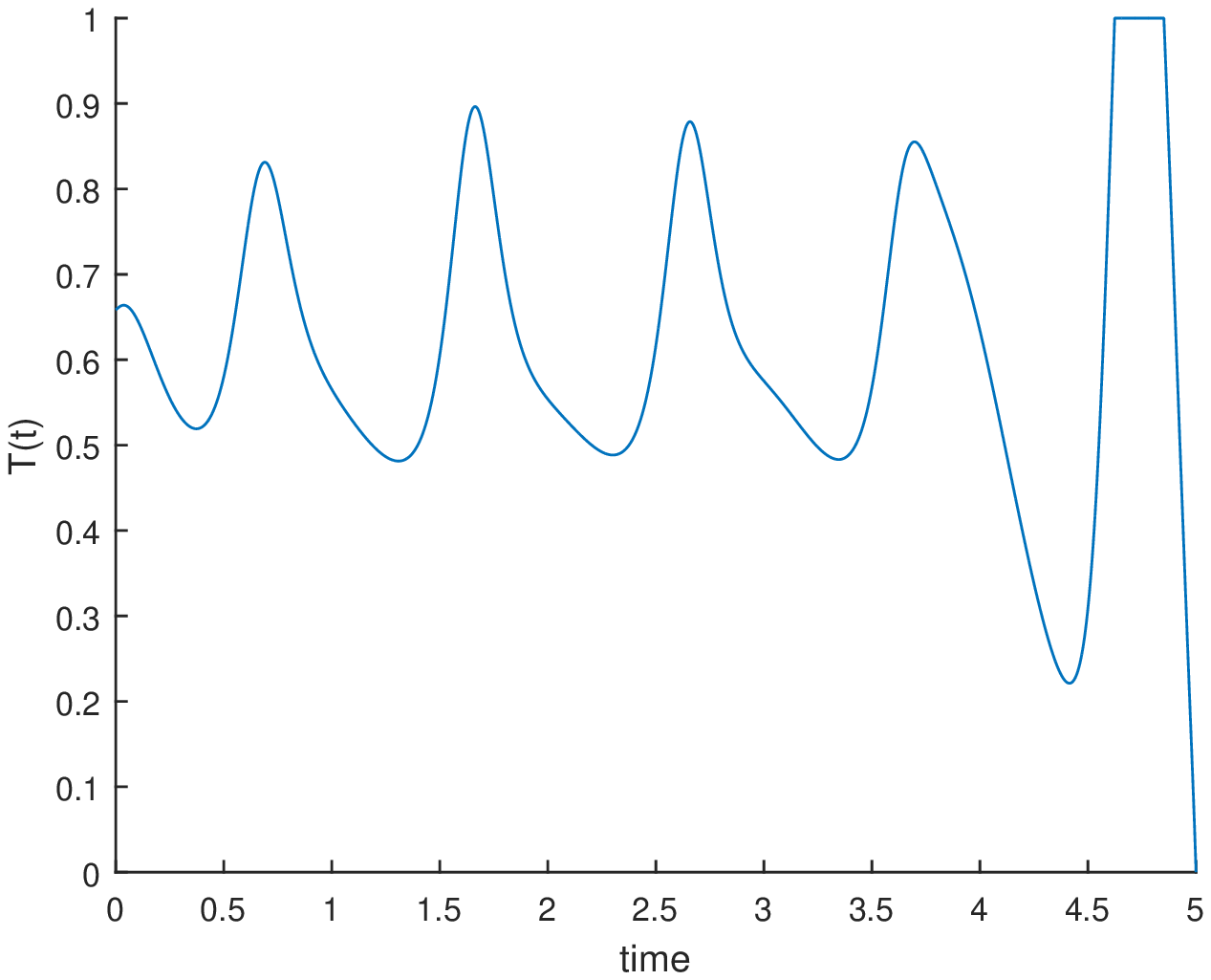}\label{fig:control_001}}\\
\subfloat[Efficacy function $F(t)$ defined in \eqref{efficacy_function}.]{
\includegraphics[scale=0.46]{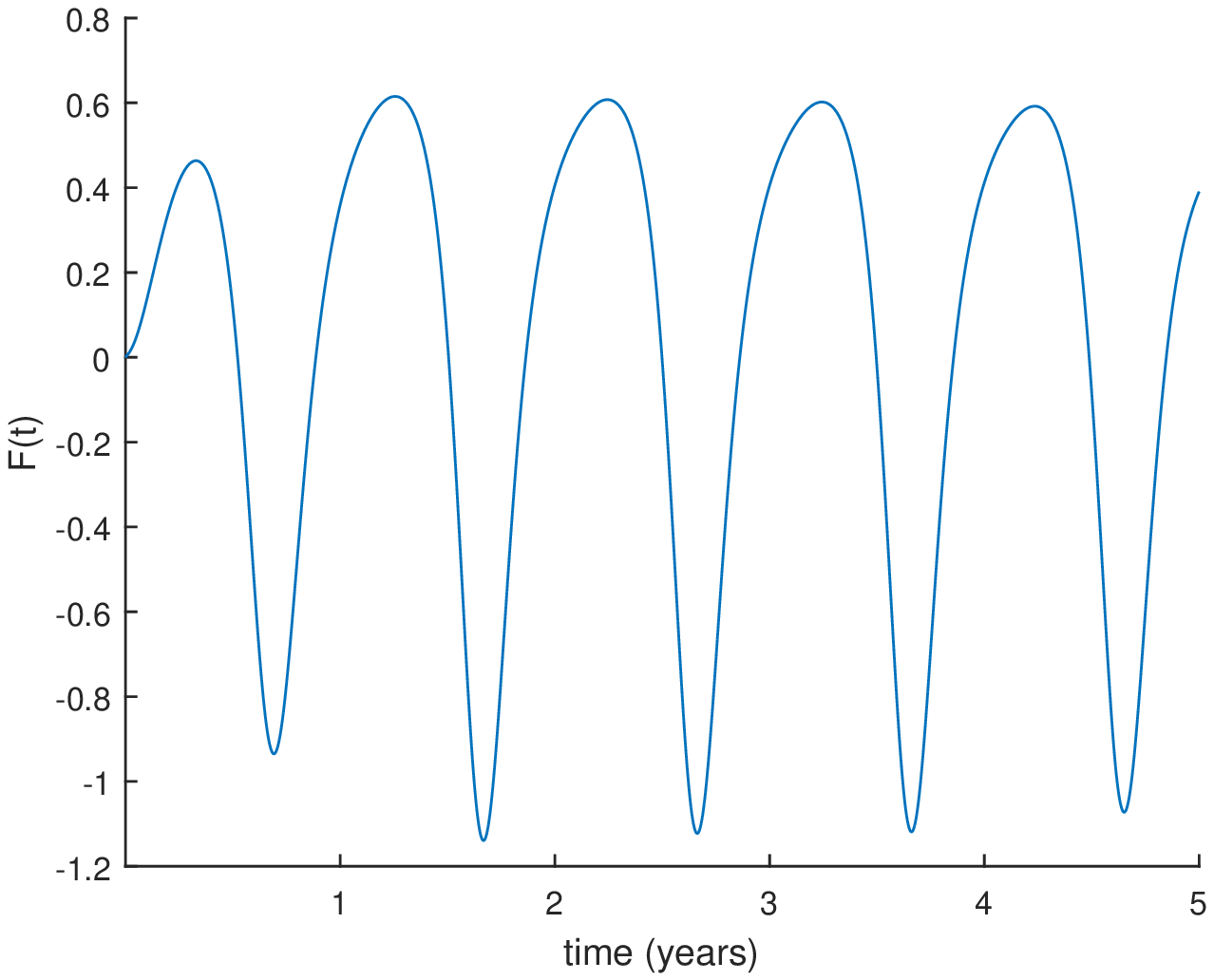}\label{fig:efficacy_001}}
\caption{Co-state variables $p_i$, $i = 1,\ldots, 4$,
extremal control \eqref{eq:ext:cont} $\mathbbm{T}$,
and efficacy function \eqref{efficacy_function} $F(t)$
associated to the fractional optimal control problem \eqref{cost-functional}--\eqref{Omega:set}
with $\alpha=0.993$ and weights $\kappa_1=1$ and $\kappa_2=0.001$.}
\label{fig:control_efficacy:K}
\end{figure}
The solution to the fractional optimal control problem is displayed in
Figures~\ref{fig:states_var:K}, \ref{fig:co_states_001} and \ref{fig:control_001}.
The periodic nature of the disease conditions the variation of the
state variables. We can also see that the control is a continuous function
with some non-regularity at the end of the time interval $[0,t_f]$. This
behaviour is motivated by the irregular oscillation of the co-state
variables, on which the control depends.

The intensity of treatment of the infectious individuals must have, periodically,
in each year of the time interval, a given period of time during which most of
the infectious individuals are treated. This ensures that the level of
infectious reach very low levels. In Figure~\ref{fig:efficacy_001}
the efficacy function \cite{rodrigues2014cost} is exhibited. It is defined by
\begin{equation}
\label{efficacy_function}
F(t)=\frac{I(0)-I^*(t)}{I(0)}=1-\frac{I^*(t)}{I(0)},
\end{equation}
where $I^*(t)$ is the optimal solution associated with the fractional optimal
control and $I(0)$ is the corresponding initial condition. This
function measures the proportional variation  in the number of
infectious individuals after the application of the control
$\mathbbm{T}^*$, by comparing the number of infected individuals
at time $t$ with the initial value $I(0)$ for which there is no
control. We observe that $F(t)$ oscillates between
$-1.14$ (lower bound) and $+0.62$ (upper bound), and  exhibits
the inverse tendency of $I(t)$.

\begin{figure}[!htb]
\centering
\subfloat[$k_1=1$ and varying $k_2$.]{
\includegraphics[scale=0.46]{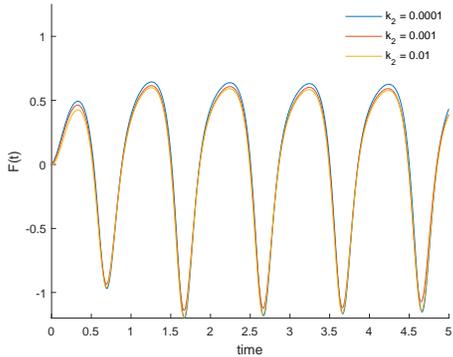}\label{fig:efficacy_vark2}}
\hspace*{1cm}
\subfloat[$k_2=0.001$ and varying $k_1$.]{
\includegraphics[scale=0.46]{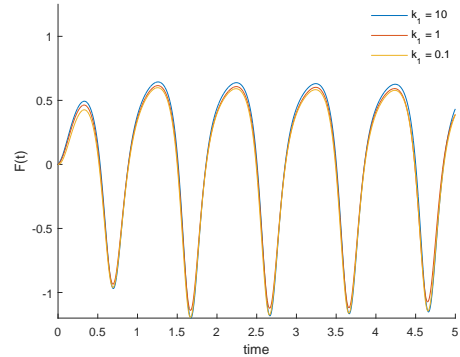}\label{fig:efficacy_vark1}}
\caption{Sensitivity analysis of the efficacy function $F(t)$ defined
by \eqref{efficacy_function} with respect to weights
$\kappa_1$ and $\kappa_2$ of the objective functional
\eqref{cost-functional}.
\emph{Left:} $\kappa_1=1$ and $\kappa_2=0.01,0.001,0.0001$.
\emph{Right:} $\kappa_2=0.001$ and $\kappa_1=0.1,1,10$.}
\label{fig:sensitivity:OF}
\end{figure}

Naturally, our results depend on the objective functional
$\mathcal{J}$, defined in \eqref{cost-functional}. They
depend, namely, on the weight constants associated with the number of
infectious individuals, $\kappa_1$, and with the cost of the treatment,
$\kappa_2$. Figure~\ref{fig:sensitivity:OF} shows that the results
do not change qualitatively by varying constants $\kappa_i$, $i=1,2$.
However, the magnitude of the efficacy changes
slightly when $\kappa_1$ and $\kappa_2$ vary independently.

To evaluate the cost and the effectiveness of the proposed fractional
control measure during the intervention period,
some summary measures are introduced.
The total cases averted by the intervention during
the time period $t_f$ is defined in \cite{rodrigues2014cost} by
\begin{equation}
\label{eq:A}
A=t_f I(0)-\int_0^{t_f}I^*(t)~dt,
\end{equation}
where $I^*(t)$ is the optimal solution associated with the fractional optimal control
$\mathbbm{T}^*$ and $I(0)$ is the corresponding initial condition. Note
that this initial condition is obtained as the equilibrium proportion
$\overline{I}$ of system \eqref{eq:modSEIRS_control} without treatment
intervention, which is independent on time, so that
$t_f I(0)=\displaystyle \int_0^{t_f}\overline{I}~dt$ represents the total infectious
cases over a given period of $t_f$ years.
Let us define effectiveness as the proportion of cases averted
on the total cases possible under no intervention \cite{rodrigues2014cost}:
\begin{equation}
\label{eq:F}
\overline{F}=\frac{A}{t_f I(0)}
=1-\frac{\displaystyle \int_0^{t_f}I^*(t)~dt}{t_f I(0)}.
\end{equation}
The total cost associated with the intervention
is defined in  \cite{rodrigues2014cost} by
\begin{equation}
\label{eq:TCI}
TC=\int_0^{t_f} C \, \mathbbm{T}^*(t)I^*(t)~dt,
\end{equation}
where $C$ corresponds to the unit cost, per person, of detection
and treatment of infectious individuals.
Following \cite{okosun2013optimal,rodrigues2014cost},
the average cost-effectiveness ratio is given by
\begin{equation}
\label{eq:ACER}
ACER=\frac{TC}{A}.
\end{equation}
\begin{table}[!htb]
\centering
\caption{Sumary of cost-effectiveness measures. Parameters according
to Tables~\ref{tab:param} and~\ref{tab:solinit}
with $\alpha=0.993$ and $C=1$.}\label{tab:efficacy}
\begin{tabular}{c@{\hspace*{1.8cm}}c@{\hspace*{1.8cm}}c@{\hspace*{1.8cm}}c}
\toprule
$A$ \eqref{eq:A} & $TC$ \eqref{eq:TCI}  & $ACER$ \eqref{eq:ACER}
& $\overline{F}$ \eqref{eq:F}\\[1mm] \midrule
172.7 & 3251.8 & 18.8 & 0.03547 \\ \bottomrule
\end{tabular}
\end{table}

The cost-effectiveness measures for the fractional optimal control problem
we have analyzed are summarized in Table~\ref{tab:efficacy}. These results
show limited effectiveness of the control treatment to reduce HRSV infectious individuals.

Another approach can be used to analyse the cost-effectiveness
of the proposed fractional optimal control problem \eqref{cost-functional}--\eqref{Omega:set}
and the classical optimal control ($\alpha=1$) investigated in \cite{rosa:delfim2018parameter},
by using the so-called \emph{incremental cost effectiveness ratio} ($ICER$) \cite{okosun2013optimal}.
This ratio is used to compare the differences between the costs and health outcomes
of two alternative intervention strategies that compete for the same resources
and it is generally described as the additional cost per additional health outcome.
First, we rank the strategies in order of increasing effectiveness, here measured
as the total infections averted $A$, defined in \eqref{eq:A}. Considering
two contending strategies $a$ and $b$, the $ICER$ of the strategy with
the least effectiveness is its ACER and for the other strategies is given by
\begin{equation}
\label{eq:ICER}
ICER(b) =\frac{TC(b)-TC(a)}{A(b)-A(a)}.
\end{equation}
\begin{table}[!htb]
\centering
\caption{Incremental cost-effectiveness ratio \eqref{eq:ICER}
for classical ($\alpha = 1$) and fractional ($\alpha = 0.993$)
HRSV disease optimal control problems. Parameters according
to Tables~\ref{tab:param} and~\ref{tab:solinit} with $C=1.$}\label{tab:icer}
\begin{tabular}{c@{\hspace*{1.7cm}}c@{\hspace*{1.7cm}}c@{\hspace*{1.7cm}}c@{\hspace*{1.7cm}}c}
\toprule
$\alpha$ & $A$ \eqref{eq:A} & $TC$ \eqref{eq:TCI}
& $ACER$ \eqref{eq:ACER} & $ICER$ \eqref{eq:ICER}\\[1mm] \midrule
1.000 & 171.1 & 3242.1 & 18.9 & 18.9\\
0.993 & 172.7 & 3251.8 & 18.8 & 6.06 \\ \bottomrule
\end{tabular}
\end{table}

Results are shown in Table~\ref{tab:icer}. The fractional order strategy
has the least $ICER$ and therefore is more cost-effective than the classical
strategy recently investigated in \cite{rosa:delfim2018parameter}.


\section{Conclusion}
\label{sec:conclu}

Human Respiratory Syncytial Virus (HRSV)
is the most common cause of lower respiratory tract infection
in infants and children worldwide.
In addition, HRSV causes serious
disease in elderly and immune compromised individuals. In this work,
we discussed fractional compartmental models for HRSV.
Estimation of the fractional order was performed for real data
of Florida from September 2011 to July 2014,
minimizing the $l_2$ norm. According to the obtained results,
the proposed models fit well the real data.
When we compare the optimal values for the fractional order SIRS-$\alpha$ and 
SEIRS-$\alpha$ models with the standard SIRS and SEIRS, 
one concludes that the absolute error of the fractional SIRS-$\alpha$ model 
reduces more than the homonym error of the fractional SEIRS-$\alpha$ model. 
Thus, we can conclude that fractional derivatives give rise to theoretical models 
that allow a significant improvement in the fitting of real data, when
compared with analogous classical models, particularly in simpler cases.
However, our results on fractional optimal control show
that treatment has a limited effect on HRSV infected individuals.
Nevertheless, a cost-effectiveness analysis of the proposed
fractional order strategy shows that it is more cost-effective
than the classical strategy followed in the literature.
As future work, we plan to investigate the usefulness of our fractional
approach in other geographical regions.


\section*{Acknowledgements}

Rosa was supported by the Portuguese Foundation
for Science and Technology (FCT) through IT (project UID/EEA/50008/2013);
Torres by FCT through CIDMA (project UID/MAT/04106/2013)
and TOCCATA (project PTDC/EEI-AUT/2933/2014
funded by FEDER and COMPETE 2020). The authors are very grateful 
to the comments and suggestions from two anonymous reviewers,
which helped them to enrich the paper. 



\end{document}